\magnification=1200
\font\titlefont=cmcsc10 at 12pt
\hyphenation{moduli}

%
%
\catcode`\@=11
\font\tenmsa=msam10
\font\sevenmsa=msam7
\font\fivemsa=msam5
\font\tenmsb=msbm10
\font\sevenmsb=msbm7
\font\fivemsb=msbm5
\newfam\msafam
\newfam\msbfam
\textfont\msafam=\tenmsa  \scriptfont\msafam=\sevenmsa
  \scriptscriptfont\msafam=\fivemsa
\textfont\msbfam=\tenmsb  \scriptfont\msbfam=\sevenmsb
  \scriptscriptfont\msbfam=\fivemsb
\def\hexnumber@#1{\ifcase#1 0\or1\or2\or3\or4\or5\or6\or7\or8\or9\or
      A\or B\or C\or D\or E\or F\fi }

\font\teneuf=eufm10
\font\seveneuf=eufm7
\font\fiveeuf=eufm5
\newfam\euffam
\textfont\euffam=\teneuf
\scriptfont\euffam=\seveneuf
\scriptscriptfont\euffam=\fiveeuf
\def\frak{\ifmmode\let\next\frak@\else
 \def\next{\errmessage{Use \string\frak\space only in math mode}}\fi\next}
\def\goth{\ifmmode\let\next\frak@\else
 \def\next{\errmessage{Use \string\goth\space only in math mode}}\fi\next}
\def\frak@#1{{\frak@@{#1}}}
\def\frak@@#1{\fam\euffam#1}
\edef\msa@{\hexnumber@\msafam}
\edef\msb@{\hexnumber@\msbfam}
\mathchardef\square="0\msa@03
\mathchardef\subsetneq="3\msb@28
\mathchardef\ltimes="2\msb@6E
\mathchardef\rtimes="2\msb@6F
\def\Bbb{\ifmmode\let\next\Bbb@\else
\def\next{\errmessage{Use \string\Bbb\space only in math mode}}\fi\next}
\def\Bbb@#1{{\Bbb@@{#1}}}
\def\Bbb@@#1{\fam\msbfam#1}
\catcode`\@=12
%
%

\def\Z{{\Bbb Z}}

\def\C{{\Bbb C}}

\def\F{{\Bbb F}}
\def\G{{\Bbb G}}
\def\V{{\Bbb V}}
\def\W{{\Bbb W}}
\def\Q{{\Bbb Q}}

\def\F{{\Bbb F}}

\def\mbar#1{{\overline{\cal M}}_{#1}}

\def\M#1{{\cal M}_{#1}}
\def\A#1{{\cal A}_{#1}}
\def\Mn#1{{\cal M}_{#1,n}}

\def\today{\ifcase\month\or
 Jan\or Febr\or  Mar\or  Apr\or May\or Jun\or  Jul\or
 Aug\or  Sep\or  Oct\or Nov\or  Dec\or\fi
 \space\number\day, \number\year}

\vskip 6.5pc
\noindent
\font\eighteenbf=cmbx10 scaled\magstep2
\vskip 2.0pc
\centerline{\eighteenbf Sur la Cohomologie des Syst\`emes Locaux}
\smallskip
\medskip
\centerline{\eighteenbf sur les Espaces des Modules des Courbes}
\smallskip
\medskip
\centerline{\eighteenbf  de Genre 2 et des Surfaces Ab\'eliennes}
\noindent
\vskip 2pc
\font\titlefont=cmcsc10 at 11pt
\centerline{\titlefont Carel Faber et Gerard van der Geer}
\vskip 2.0pc
\centerline{\bf On the cohomology of local systems on the moduli }
\centerline{\bf spaces of curves of genus $2$ and of abelian surfaces}
\smallskip
\centerline{\bf Summary}
We consider the cohomology of local systems on the moduli space
of curves of genus 2 and the moduli space of abelian surfaces.
We  give an explicit formula for the Eisenstein cohomology and a conjectural
formula for the endoscopic contribution. We show how counting curves
over finite fields provides us with detailed information about
Siegel modular forms.
\medskip
\centerline{\bf \S0 R\'esum\'e}
\smallskip
L'objet est la cohomologie des syst\`emes locaux sur les espaces
$\M2$ des modules des courbes de genre $2$ et $\A2$ des modules des surfaces
ab\'eliennes. Nous donnons une formule explicite pour la
cohomologie d'Eisenstein et une formule conjecturale pour la
contribution endoscopique. Notre calcul des courbes sur des corps finis
donne des renseignements pr\'ecis sur les formes modulaires de Siegel.
\medskip
\centerline{\bf \S1 Syst\`emes Locaux et Caract\'eristiques Motiviques
d'Euler}
\smallskip
Soit $\M2$ l'espace des modules des courbes de genre $2$ et soit
$\A2$ l'espace des modules des surfaces ab\'eliennes principalement
polaris\'ees.
Ce sont des champs alg\'ebriques d\'efinis sur $\Z$. En associant
\`a une courbe sa jacobienne on obtient une immersion ouverte  $\M2 \to
\A2$. Nous notons ${\cal A}_{1,1}=\A2-\M2$.
La courbe universelle $\pi': {\cal C} \to \M2$ et la surface
ab\'elienne universelle $\pi:{\cal X}\to \A2$ d\'efinissent des syst\`emes
locaux $\V'=R^1\pi'_*(\Q)$ sur $\M2$ et $\V =R^1\pi_*(\Q)$ sur $\A2$
et sous l'immersion $j:\M2 \to \A2$ on a $j^*\V=\V^{\prime}$.
On a un accouplement
symplectique $\V\otimes \V \to \Q(-1)$.

A chaque paire d'entiers $(l,m)$ avec $l\geq m\geq 0$ on peut
associer une repr\'esentation irr\'eductible du groupe ${\rm Sp}(4,\Q)$
et on la rel\`eve en une repr\'esentation de ${\rm GSp}(4,\Q)$ de poids dominant
$(\ell -m)\gamma_{\beta}+m \gamma_{\alpha} -(\ell+m) \eta$ avec
$\gamma_{\alpha}=(1,1)$
et $\gamma_{\beta}=(0,1)$ deux racines fundamentales et $\eta$
le multiplicateur. Ainsi $(\ell,m)=(1,0)$ donne le contragr\'edient
de la repr\'esentation standard. \c{C}a d\'efinit
un syst\`eme local $\V_{l,m}$ sur $\A2$,
facteur direct de ${\rm Sym}^{\ell-m}(\V)\otimes
{\rm Sym}^m(\wedge^{2}(\V))$ de poids $\ell+m$. On a $\V= \V_{1,0}
$ et $\wedge^2 \V= \V_{1,1}\oplus \V_{0,0}(-1)$.
Une paire $(l,m)$ est appel\'ee r\'eguli\`ere si on a $l>m>0$.

On s'int\'eresse \`a la caract\'eristique motivique d'Euler
$e_c(\A2,\V_{l,m})$ d\'efinie par
$$
e_c(\A2,\V_{l,m})= \sum_{i=0}^6 (-1)^i [H^i_c(\A2,\V_{l,m})],
$$
o\`u on prend la classe de la cohomologie \`a support compact dans $K_0$
d'une cat\'egorie convenable, par exemple, la cat\'egorie
des modules de Hodge mixtes sur $\A2$ ou des
motifs effectifs de Chow.
Pour $l+m$ impair  on voit que $e_c(\A2,\V_{l,m})=0$,
car $-{\rm Id}\in{\rm Sp}(4)$ agit comme $-1$ sur $\V$.
Dans le cas de genre $1$ on a
\'egalement l'espace $\A1$ des modules
des courbes elliptiques, sa courbe universelle
$\pi : {\cal E} \to \A1$, et le syst\`eme local $\W=R^1\pi_*(\Q)$. On pose
$\W_k={\rm Sym}^k(\W)$. Dans ce cas  on sait que
$e_c(\A1,\W_k)=0$
pour $k$ impair,
et $e_c(\A1,\W_k)= -S[k+2] -1$, pour $k\geq 2$ pair, o\`u
$S[k+2]$ est le motif des formes paraboliques
de poids $k+2$ sur ${\rm SL}(2,\Z)$, cf.\ [D], [Sch].
On a $e_c(\A1,\W_0)=L$ o\`u $L$ est le motif de Tate de poids $2$.
\medskip
\centerline{\bf \S2 La Contribution de la Cohomologie d'Eisenstein}
\smallskip
On a une fl\`eche naturelle $H_c^*(\A2, \V_{l,m}) \to H^*(\A2,\V_{l,m})$.
L'image est appel\'ee la cohomologie int\'erieure $H_{!}^*$ et on note
sa caract\'eristique motivique d'Euler $e_!$.
On d\'efinit la cohomologie d'Eisenstein
comme le noyau
et la caract\'eristique correspondante est d\'efinie par
$e_{\rm Eis}(\A2,\V_{l,m})=e_c(\A2,\V_{l,m})-e_{!}(\A2,\V_{l,m})$.
En utilisant les techniques et r\'esultats de Harder [H1,2], Pink [P]
et Schwermer [Schw]
on peut d\'emontrer le th\'eor\`eme suivant pour la cohomologie de Betti
et pour la cohomologie \'etale avec coefficients  $\ell$-adiques
en caract\'eristique $p>0$.
\proclaim Th\'eor\`eme. Soit $(l,m)$ r\'eguli\`ere.
La caract\'eristique motivique d'Euler
pour la cohomologie de Betti et pour la
cohomologie \'etale $\ell$-adique de la cohomologie d'Eisen\-stein est
$$
e_{\rm Eis}(\A2,\V_{l,m})=
s_{l-m+2}-s_{l+m+4}L^{m+1}
+\cases{
S[m+2]+1 &\quad $l$ pair, \cr -S[l+3] &\quad $l$ impair,}
$$
o\`u $s_n$ est la dimension de l'espace vectoriel des formes paraboliques
de poids $n$ sur ${\rm SL}(2,\Z)$.
\par
Pour $l>m=0$ et pour $l=m>0$
cette formule garde un sens quand on
interpr\`ete $S[2]$ comme $-L-1$ et $s_2$ comme $-1$ et devrait \^etre vraie
sauf dans le cas $l=m>0$, $m$ pair quand il y ont des annullements
inattendus des fonctions $L$  des formes
modulaires, cf.\ [H2].
\medskip
\centerline{\bf \S3 Formes Modulaires de Siegel}
\smallskip
Soit $(U,\rho)$ une repr\'esentation irr\'eductible de ${\rm GL}(2,\C)$
de type $(j,k)$, c'est-\`a-dire de la forme
${\rm Sym}^j(R)\otimes \det (R)^k$,
o\`u $R$ est la repr\'esentation standard.
Une forme modulaire de Siegel de poids
$\rho$ ou $(j,k)$ est une fonction holomorphe $f: {\cal H}_2 \to U$
sur le demi-plan de Siegel ${\cal H}_2$ telle que
$$
f((A\tau +B)(C\tau +D)^{-1})=\rho(C\tau +D) f(\tau)
$$
pour $\tau \in {\cal H}_2$ et $(A,B; C, D) \in {\rm Sp}(4,\Z)$, cf.\ [A].
Notons $M_{j,k}$ l'espace vectoriel de telles formes de poids $(j,k)$ et
$S_{j,k}$ le sous-espace des formes paraboliques.
On a $M_{j,k}=(0)$ si $j$ est impair, ou $j<0$ ou $k<0$.
Les formes modulaires forment un anneau
${\cal R}=\oplus_{j,k} M_{j,k}$ et $\oplus S_{j,k}$ est un id\'eal de
cet anneau.
Le sous-anneau ${\cal R}^0=\oplus_{k} M_{0,k}$ est la $\C$-alg\`ebre des formes
modulaires classiques et des g\'en\'erateurs de ${\cal R}^0$ ont \'et\'e
d\'etermin\'es par Igusa.
On connait \'egalement des g\'en\'erateurs du
${\cal R}^0$-modules $\oplus_{k} M_{2,k}$ et $\oplus_{k} M_{4,k}$
d'apr\`es Satoh [S] et Ibukiyama [I1,I2].
Mais \`a l'exception de ces r\'esultats
et une formule pour la dimension
de $S_{j,k}$ due \`a Tsushima [Ts],  presque rien n'est connu.

L'espace $S_{6,8}$ est de dimension $1$ et \`a notre demande
Ibukiyama [I3] a construit
une forme $0\neq F \in S_{6,8}\,$, en utilisant le r\'eseau
$\Gamma =\{ x \in \Q^{16}\, : \, 2x_i\in \Z, x_i-x_j \in \Z,
\sum_{i=1}^{16} x_i \in 2\Z\} $.
On pose $a=(2,i,i,i,i,0,\ldots,0) \in \C^{16}$ et on note par $(\, , \, )$
le produit scalair usuel. Soit $F=(F_0,\ldots,F_6)$ le vecteur de fonctions
sur ${\cal H}_2$ d\'efini par
$$
F_{\nu}=\sum_{x,y \in \Gamma} (x,a)^{6-\nu} (y,a)^{\nu} e^{\pi i ((
x,x)\tau_1+(x,y)\tau_2+(y,y)\tau_3)} \qquad (\nu=0,\ldots,6)
$$
pour $\tau= (\tau_1, \tau_2; \tau_2, \tau_3) \in {\cal H}_2$.
Le r\'esultat de Ibukiyama dit que $F\neq 0$ et $F \in S_{6,8}$.
\medskip
\centerline{\bf \S4 Syst\`emes Locaux et Formes Modulaires}
\smallskip
On sait que pour $(l,m)$ r\'eguli\`ere $H^i_{!}(\A2,\V_{l,m})=(0)$
quand $i\neq 3$, cf.\ [T].
Faltings a montr\'e dans [F] (cf.\ [F-C], VI, Th.\ 5.5) que
$H^3(\A2, \V_{l,m})$ et $H^3_c(\A2, \V_{l,m})$ sont munis de filtrations
de Hodge de poids $\geq l+m+3$ et $\leq l+m+3$ resp.\ et on trouve
que $H_{!}^3(\A2, \V_{l,m})$  est muni d'une filtration de Hodge
$$
(0) \subset F^{l+m+3} \subset F^{l+2} \subset F^{m+1} \subset F^{0} =
H^3_{!}(\A2, \V_{l,m}).
$$
De plus, on sait que $F^{l+m+3} \cong S_{l-m,m+3}$.
Soit $M$ le groupe `endoscopique' ${\rm GL}(2)\times
{\rm GL }(2)/\G_m $, o\`u $\G_m\to {\rm GL}(2)\times{\rm GL }(2)$ via
$x \mapsto (x,x^{-1})$.
Le rel\`evement endoscopique de $M$ contribue aussi
\`a la cohomologie int\'erieure $H^3_{!}(\A2,\V_{l,m})$.
Nous n'avons pas pu trouver une assertion pr\'ecise sur cette
contribution dans la litt\'erature. Nous nous attendons \` a ce
que les experts sachent
d\'emontrer notre description conjecturale de cette contribution
endoscopique. Cette description semble en accord avec les r\'esultats
de Kudla et Rallis [K-R].

\proclaim Conjecture. Si la paire $(l,m)$ est reguli\`ere la
contribution endoscopique du groupe  $M$
\`a $e_!(\A2,\V_{l,m})$
a une intersection nulle
avec $F^{l+m+3}$ et est \'egale \`a
$$
e_{\rm endo }(\A2,\V_{l,m})= -s_{l+m+4} S[l-m+2] L^{m+1}.
$$
\par
Comme nous le verrons nos calculs donnent
un tr\`es fort  support en faveur de  cette conjecture.
Nous d\'efinissons
$$
S[l-m,m+3]:= H_!^3(\A2, \V_{l,m}) - H^3_{\rm endo}(\A2, \V_{l,m}),
$$
o\`u $H^3_{\rm endo}$ est la partie endoscopique de $H^3_{!}$.
Cela devrait \^etre un motif de rang \'egal \`a
$4\dim S_{l-m,m+3}$ avec des poids
de Hodge $l+m+3$, $l+2$, $m+1$ et $0$. Nous le consid\'erons
dans la cat\'egorie
des structures de Hodge mixtes ou des repr\'esentations galoisiennes.
On a la formule conjecturale
$$
e_c(\A2,\V_{l,m})=
-S[l-m,m+3]-s_{l+m+4} S[l-m+2] L^{m+1}+e_{\rm Eis}(\A2,\V_{l,m}).
$$
Cette formule garde un sens pour $l>m=0$ et pour $l=m>0$ quand on
interpr\`ete $S[2]$ comme $-L-1$ et $s_2$ comme $-1$.
\medskip
\centerline{\bf \S5 Calcul des Points sur des Corps Finis}
\smallskip
Soit $p$ un nombre premier et $q=p^e$. Comme $\A2$ et $\M2$ sont d\'efinis sur
$\Z$ on peut consid\'erer $\A2\otimes \F_p$ et $\M2\otimes \F_p$
et on peut d\'efinir la trace $t(l,m,q)$ de la $e$-i\`eme
puissance de Frobenius pour la variante $\ell$-adique de $\V_{l,m}$
sur $\A2\otimes \F_p$, cf.~[Be].

Pour certains $q$ et une famille convenable
de courbes ${\cal D} \to N$ de genre $2$ avec $N \to \M2$
fini nous avons d\'etermin\'e les points rationnels $N(\F_q)$
et les polyn\^omes de Weil des courbes correspondantes.
Ainsi nous avons calcul\'e les fr\'equences des polyn\^omes de Weil.
(Cf.\ le cas $g=1$, o\`u Birch [Bi] a fait des calculs analogues.)
Divisant par le degr\'e du recouvrement $N\to \M2$ on trouve
les fr\'equences pour le champ alg\'ebrique
$\M2$.
Ce sont des nombres rationnels.
(De fa\c{c}on equivalente: prendre un
repr\'esentant $X$ de chaque classe
d'isomorphie sur $\F_q$, et sommer les
${\rm Tr}(F,V_{\ell,m}(H^1))/\#{\rm Aut}_{\F_q}(X)$.
Notons qu'une courbe sur un corps fini
peut \^etre d\'efinie sur son corps de modules.)
En y ajoutant la contribution
des courbes stables de genre $2$ avec deux composantes elliptiques
on peut d\'eterminer les traces $t(l,m,q)$ pour ces $q$. Nous l'avons
fait pour $q \leq 37$, $q\neq 27$.

Supposant ainsi notre conjecture vraie
nous pouvons calculer les traces des op\'erateurs de Hecke $T(p)$
sur les espaces $S_{j,k}$ pour tous les nombres premiers
$p \leq 37$.
Notons qu'en g\'en\'eral, m\^eme si on a une forme explicite
comme pour $F \in S_{6,8}$, il est tr\`es difficile de calculer les
coefficients de Fourier et les valeurs propres des op\'erateurs de Hecke.
\medskip
\centerline{\bf \S6 Polyn\^omes Caract\'eristiques de la Forme Modulaire $F$}
\smallskip
Ici et dans la section suivante nous supposons vraie la conjecture sur la
partie endoscopique.
Pour une forme propre des op\'erateurs de Hecke $T(n)$ dans $S_{l-m,m+3}$
avec valeurs propres $\lambda(n)$,
le polyn\^ome caract\'eristique de Frobenius est donn\'e
par
$$
1-\lambda(p)X+(\lambda(p)^2-\lambda(p^2)-p^{l+m+2})X^2
-\lambda(p)p^{l+m+3}X^3+ p^{2l+2m+6}X^4,
$$
cf.~[A], p.~164.
Nous donnons les polyn\^omes caract\'eristiques  de Frobenius pour $p$ petit
pour la forme $F\in S_{6,8}$ et les pentes de leurs polygones de
Newton.
En utilisant les calculs des coefficients de Fourier de $F$ de Ibukiyama,
on peut controler $\lambda(2)=0$.
\smallskip
\vbox{
\bigskip\centerline{\def\quad{\hskip 0.6em\relax}
\def\quod{\hskip 0.5em\relax }
\vbox{\offinterlineskip
\hrule
\halign{&\vrule#&\strut\quod\hfil#\quad\cr
height2pt&\omit&&\omit&&\omit&&\omit&\cr
&$p$&&$\lambda(p)$&&$\lambda(p^2)$&&{\rm pentes}&\cr
height2pt&\omit&&\omit&&\omit&&\omit&\cr
\noalign{\hrule}
height2pt&\omit&&\omit&&\omit&&\omit&\cr
&$2$&&$0$&&$-57344$&&$13/2,25/2$&\cr
&$3$&&$-27000$&&$143765361$&&$3,7,12,16$&\cr
&$5$&&$2843100$&&$-7734928874375$&&$2,7,12,17$&\cr
&$7$&&$-107822000$&&$4057621173384801$&&$0,6,13,19$&\cr
height2pt&\omit&&\omit&&\omit&&\omit&\cr
} \hrule}
}}
\medskip
\centerline{\bf \S7 Exemples de Valeurs Propres des Op\'erateurs de Hecke}
\smallskip
Il y a 29 cas de paires r\'eguli\`eres $(l,m)$ tel que $\dim S_{l-m,m+3}=1$.
Dans ces cas les traces de $T(p)$ donnent les valeurs propres et nos
r\'esultats permettent de calculer ces valeurs propres pour
$p\leq 37$.
Pour le cas $(l,m)=(13,11)$ nos r\'esultats sont en accord avec les
calculs de Satoh de quelques valeurs propres sur $S_{2,14}$, cf.\ [S], p.~351.
Pour illustrer ces r\'esultats
nous donnons les valeurs propres $\lambda(p)$ des op\'erateurs de Hecke $T(p)$
sur les espaces $1$-dimensionnels $S_{8,8}$ et $S_{12,6}$ pour $p\leq19$.
\smallskip
\vbox{
\bigskip\centerline{\def\quad{\hskip 0.6em\relax}
\def\quod{\hskip 0.5em\relax }
\vbox{\offinterlineskip
\hrule
\halign{&\vrule#&\strut\quod\hfil#\quad\cr
height2pt&\omit&&\omit&&\omit&\cr
&$p$&& $\lambda(p)$  sur $S_{8,8}$&&$\lambda(p)$ sur $S_{12,6}$&\cr
height2pt&\omit&&\omit&&\omit&\cr
\noalign{\hrule}
height2pt&\omit&&\omit&&\omit&\cr
&$2$&&$2^6\cdot 3 \cdot 7$&&$-2^4\cdot 3 \cdot 5$&\cr
&$3$&&$-2^3\cdot 3^2 \cdot 89$&&$2^3 \cdot 3^5 \cdot 5 \cdot 7$&\cr
&$5$&&$-2^2\cdot 3 \cdot 5^2 \cdot 13^2 \cdot 607$&&$2^2\cdot 3 \cdot 5^2 \cdot 7 \cdot 79 \cdot 89$&\cr
&$7$&&$2^4\cdot 7 \cdot 109 \cdot 36973$&&
 $-2^4 \cdot 5^2 \cdot 7 \cdot 119633$&\cr
&$11$&&$2^3\cdot 3 \cdot 4759 \cdot 114089$&&
$ 2^3 \cdot 3 \cdot 23 \cdot 2267 \cdot 2861$&\cr
&$13$&&$-2^2\cdot 13 \cdot 17 \cdot 109 \cdot 3404113$&&
$2^2 \cdot 5 \cdot 7 \cdot 13 \cdot 50083049$&\cr
&$17$&&$2^2\cdot 3^2 \cdot 17 \cdot 41 \cdot 1307 \cdot 168331$&&$
-2^2\cdot 3^2 \cdot 5 \cdot 7 \cdot 13 \cdot 47 \cdot 14320807$&\cr
&$19$&&$-2^3\cdot 5 \cdot 74707 \cdot 9443867$&&
$-2^3 \cdot 5 \cdot 7^3 \cdot 19 \cdot 2377 \cdot 35603$&\cr
height2pt&\omit&&\omit&&\omit&\cr
} \hrule}
}}
\medskip
\centerline{\bf \S8 Les Caract\'eristiques Motiviques d'Euler de $\Mn2$}
\smallskip
Soit $\Mn1$ (resp.\ $\Mn2$)
l'espace des modules des courbes lisses de genre $1$ (resp.\ $2$) avec $n>0$
(resp.\ $n\geq 0$) points distincts ordonn\'es. C'est un champ
alg\'ebrique de dimension $n$ (resp.\ $3+n$) sur $\Z$.
Soit $\mbar{2,n}$ la compactification qui donne les modules
des courbes stables
avec $n$ points. On consid\`ere les caract\'eristiques motiviques d'Euler
$e_c(\Mn2)$ et $e_c(\mbar{2,n})$.
Getzler [G1] a d\'etermin\'e ces caract\'eristiques pour $n\leq3$.

Il y a une relation entre
$e_c(\Mn1)$ et les $e_c(\A1,\W_k)$ et Getzler [G2] l'a
\'ecrite explicitement sous la forme \'el\'egante
$$
e_c({\cal M}_{1,n+1})/n!=
{\rm Res}_0\big[ {L-\omega -L/\omega \choose n } \big(
\sum_{k=1}^{\infty}
({S[2k+2]+1 \over L^{2k+1}} )\omega^{2k} -1\big)
\cdot (\omega - L/\omega) d\omega \big]
$$
o\`u $L$ est le motif de Tate et ${\rm Res}_0$ veut dire le residu en $
\omega=0$.
Pour $g=2$ on a une relation analogue, aussi due a Getzler et cela
nous permet de calculer
$e_c(\Mn2)$ quand on sait les $e_c(\A2,\V_{k,l})$
et les $e_c({\cal A}_{1,1},\V_{k,l})$ avec $k+l\leq n$.

Nous avons calcul\'e pour $p\leq 233$ et $q=2^k$ avec $1\leq k \leq 7$
les fr\'equences des courbes avec un nombre donn\'e de
points rationels sur $\F_p$
et $\F_q$ et cela nous a permis de calculer $\# \Mn2(\F_p)$
et $\# \Mn2(\F_q)$. Cela fournit un contr\^ole pour les
$e_c(\Mn2)$
que nous avons trouv\'es.
Pour $n\leq 16$ nous connaissons la partie endoscopique
et par consequent nous connaissons
$e_c(\Mn2)$.
Nous donnons quelques exemples.
$$
\eqalign{
e_c({\cal M}_{2,10})= &
 L^{13}+ 10L^{12}-120L^{11}-420L^{10}+8253L^9
                -38931L^8+95927L^7+\cr
&\quad -156313L^6+  212730L^5 -189334L^4-166663L^3+604236L^2+\cr
&\qquad \qquad \hfill  -233280 L-302400+ (L-9)S[12].\cr
e_c({\cal M}_{2,16}) = & A(L)+B(L)S[18]+C(L)S[16]+D(L)S[12]+
2548 \, e_c({\cal M}_2 ,{\V }_{11,5})\cr
}
$$
o\`u
$$
e_c({\cal M}_{2}, {\V}_{11,5})=-S[6,8]-(L+1)S[12]-L^5-2L^4-2L^3-2L^2-2L-2
$$
avec $S[6,8]$ le `motif' des formes modulaires de Siegel de poids $(6,8)$ et
$$
\eqalign{
A(L)=L^{19}+16L^{18}-560L^{17}-9828L^{16}+441714L^{15}-6084744L^{14}+
41896790L^{13}+&\cr -109303012L^{12}-635722737L^{11}+8143511948L^{10}
-43836908908L^9+&\cr 144755642044 L^8-279604866542L^7+109307474312L^6+ &\cr
1000400749388L^5-2508087212954L^4+ 1168650933424L^3+&\cr
3250035688136L^2 -2584542638104L-1743565818904&\cr
}
$$
et
$$
\eqalign{
B(L)&=-51480L^2+ 649792L-2102115\cr
C(L)&=-1560L^2+19175L -64260\cr
D(L)&=8008L^7-264264L^6+4086368L^5-39326716L^4+ 256866324L^3+\cr
&\qquad -1125323276L^2+ 2979292862L-3563125592.\cr
}
$$
Vu les grands coefficients de ces formules il est pr\'ef\'erable
de travailler avec les $e_c({\cal A}_2, \V_{l,m})$ et
les $e_c({\cal M}_2,\V_{l,m})$.
\bigskip
Nous remercions P.\ Deligne, G.\ Harder, A.J.\ de Jong et J.-P.\ Serre
de leurs remarques sur une version pr\'eliminaire de ce manuscrit
et T.\ Ibukiyama pour correspondance utile.
Enfin, nous remercions
S.\ del Ba{\~n}o, notre collaborateur dans le premier phase de ce projet.
\medskip
\centerline{\bf Bibliographie}
\smallskip
\noindent
[A] T.\ Arakawa: Vector valued Siegel's modular forms of degree $2$
and the associated Andrianov $L$-functions.
{\sl Manuscr.\ Math.\ \bf 44} (1983), 155--185.
\smallskip
\noindent
[Be] K.\ Behrend: The Lefschetz trace formula for algebraic stacks.
{\sl Invent.\ Math.\ {\bf 112} }(1993), no. 1, 127--149.
\smallskip
\noindent
[Bi] B.\ Birch: How the number of points of an elliptic curve over a
fixed prime field varies. {\sl J.\  London Math.\  Soc.\ \bf  43} (1968),
   57--60.
\smallskip
\noindent
[D] P.\ Deligne: Formes modulaires et repr\'esentations $\ell$-adiques. S\'em.\ Bourbaki exp.\ 355, Lecture Notes in Math.\ 179 (vol.\ 1968/1969).
Springer Verlag, Berlin, 1971.
\smallskip
\noindent
[F] G.\ Faltings: On the cohomology of locally symmetric Hermitian spaces.
Paul Dubreil and Marie-Paule Malliavin algebra seminar, Paris, 1982, 55--98,
Lecture Notes in Math., 1029,
Springer, Berlin, 1983.
\smallskip
\noindent
[F-C] G.\ Faltings, C.L.\ Chai: Degeneration of abelian varieties.
Ergebnisse der Math. 22. Springer Verlag.
\smallskip
\noindent
[G1] E.\ Getzler: Topological recursion relations in genus $2$.
Integrable systems and algebraic geometry (Kobe/Kyoto, 1997),
   73--106, World Sci.~Publishing, River Edge, NJ, 1998.
\smallskip
\noindent
[G2] E.\ Getzler: Resolving mixed Hodge modules on configuration spaces.
{\sl Duke Math.\ J.\ \bf 96} (1999),  175--203.
\smallskip
\noindent
[H1] G.\ Harder:
 Eisensteinkohomologie und die Konstruktion gemischter Motive.
Lecture Notes in Mathematics, 1562. Springer-Verlag, Berlin, 1993.
\smallskip
\noindent
[H2] G.\ Harder: Modular symbols and special values of automorphic
$L$-functions. Manu\-scrit non-publi\'e.
\smallskip
\noindent
[I1] T.\ Ibukiyama: Vector valued Siegel modular forms of symmetric
tensor representations of degree $2$. Manuscrit non-publi\'e.
\smallskip
\noindent
[I2] T.\ Ibukiyama: Vector valued Siegel modular forms
of ${\rm det}^k{\rm Sym}(4)$ and ${\rm det}^k
{\rm Sym}(6)$. Manuscrit non-publi\'e.
\smallskip
\noindent
[I3] T.\ Ibukiyama: Lettre \`a G.\ van der Geer, Juillet 2001.
\smallskip
\noindent
[K-R] S.\ Kudla, S.\ Rallis: A regularized Siegel-Weil formula: The first term identity. {\sl Ann.\ of Math.\ \bf 140} (1994), 1--80.
\smallskip
\noindent
[P] R.\ Pink: On $l$-adic sheaves on Shimura varieties and their
higher direct images in the Baily-Borel compactification. {\sl Math.\ Ann.\
\bf 292} (1992),  197--240.
\smallskip
\noindent
[S] T.\  Satoh: On certain vector valued Siegel
modular forms of degree 2. {\sl Math.\ Ann.\
\bf 274} (1986), 335--352.
\smallskip
\noindent
[Sch] A.\ Scholl: Motives for modular forms. {\sl Invent.\  Math.\ \bf 100}
 (1990), no. 2, 419--430.
\smallskip
\noindent
[Schw] J.\ Schwermer: On Euler products and residual Eisenstein
cohomology classes for Siegel modular varieties.
{\sl Forum Math.\ \bf 7} (1995), 1--28.
\smallskip
\noindent
[T] R.\ Taylor: On the $\ell$-adic cohomology of Siegel threefolds.
{\sl Invent.\ Math.\ \bf 114} (1993), 289--310.
\smallskip
\noindent
[Ts] R.\ Tsushima: An explicit dimension formula for the spaces of
generalized automorphic forms with respect to ${\rm Sp}(2,\Z)$.
{\sl Proc.\ Jap.\ Acad.\ \bf 59A} (1983), 139--142.

\bigskip
\noindent Institutionen f\"or Matematik \par
\noindent Kungl Tekniska H\"ogskolan \par
\noindent S-100 44 Stockholm, Su\`ede \par
\noindent faber@math.kth.se \par

\bigskip
\noindent Korteweg-de Vries Instituut \par
\noindent Universiteit van Amsterdam \par
\noindent Plantage Muidergracht 24 \par
\noindent NL-1018 TV Amsterdam, Pays-Bas \par
\noindent geer@science.uva.nl \par
\end